\documentclass[amsfonts]{article}
\usepackage{amssymb}

\newtheorem{theorem}{Theorem}[section]

\newtheorem{lemma}[theorem]{Lemma}

\def\bR{\mathbb{R}}
\def\bZ{\mathbb{Z}}

\def\cB{\mathcal{B}}

\def\cH{\mathcal{H}}

\def\cP{\mathcal{P}}

\begin{document}

\title{A Note on a Fenyman-Kac-Type Formula}

\author{Raluca M. Balan\thanks{Research supported by a grant from
the Natural Sciences and Engineering Research Council of Canada.} }

\date{April 28, 2009}

\maketitle

\begin{abstract}
\noindent In this article, we establish a probabilistic
representation for the second-order moment of the solution of
stochastic heat equation in $[0,1] \times \bR^d$, with
multiplicative noise, which is fractional in time and colored in
space. This representation is similar to the one given in
\cite{DMT08} in the case of an s.p.d.e. driven by a Gaussian noise,
which is white in time. Unlike the formula of \cite{DMT08}, which is
based on the usual Poisson process, our representation is based on
the planar Poisson process, due to the fractional component of the
noise.
\end{abstract}




\section{Introduction}

The classical Feynman-Kac (F-K) formula gives a stochastic
representation for the solution of the heat equation with potential,
as an exponential moment of a functional of Brownian paths (see e.g.
\cite{karatzas-shreve91}). This representation is a useful tool in
stochastic analysis, in particular for the study stochastic partial
differential equations (s.p.d.e.'s). We mention briefly several
examples in this direction. The F-K formula lied at the origin of
the existence, large time asymptotic and intermittency results of
\cite{carmona-molchanov94} and \cite{carmona-viens98}, for the
solution of the heat equation with random potential on $\bZ^d$,
respectively $\bR^d$. The same method, based on a discretized F-K
formula, was used in \cite{tindel-viens02} for obtaining an upper
bound for the exponential behavior of the solution of the heat
equation with random potential on a smooth compact manifold. The
technique of \cite{carmona-molchanov94} was further refined in
\cite{CMS05} in the case of parabolic equations with L\'evy noise,
for proving the exponential growth of the solution. The F-K formula
was used in \cite{DKPW02} for solving stochastic parabolic
equations, in the context of white noise analysis. A F-K formula for
the solution of the stochastic KPP equation is used in \cite{OVZ00}
for examining the asymptotic behavior of the solution.

The present work has been motivated by the recent article
\cite{DMT08}, in which the authors obtained an alternative
probabilistic representation for the  solution of a deterministic
p.d.e., as well as a representation for the moments of the
(mild-sense) solution of a s.p.d.e. perturbed by a Gaussian noise
$\dot F$, with ``formal'' covariance:
$$E[\dot{F}_{t,x}\dot{F}_{s,y}]=\delta_0(t-s)f(x-y).$$
More precisely, in \cite{DMT08}, $\{F(h), h \in \cP\}$ is a
zero-mean Gaussian process with covariance $E(F(h)F(g))=\langle h,g
\rangle_{\cP}$, where $\cP$ is the completion of $\{1_{[0,t] \times
A}; t \in [0,1], A \in \cB_b(\bR^d)\}$ with respect to the inner
product $\langle \cdot,\cdot \rangle_{\cP}$ given by
$$\langle
1_{[0,t] \times A},1_{[0,s] \times B} \rangle_{\cP}=(t \wedge s)
\int_{A} \int_{B} f(x-y)dydx.$$ (Here, $\cB_{b}(\bR^d)$ denotes the
class of bounded Borel sets in $\bR^d$.)

In the particular case of the stochastic heat equation:
\begin{eqnarray}
\label{heat-Dalang}
\frac{\partial u}{\partial t} &=& \frac{1}{2} \Delta u + u \ \dot {F} , \quad t  >0, x\in \bR^{d} \\
\nonumber u_{0,x} &=& u_0(x), \quad x \in \bR^d,
\end{eqnarray}
the representation for the second-moments of the (mild-sense)
solution $u$ is:
\begin{equation}
\label{Dalang-representation} E[u_{t,x}u_{t,y}]=e^t
E_{x,y}\left[w(t-\tau_{N_t},B_{\tau_{N_t}}^{1})
w(t-\tau_{N_t},B_{\tau_{N_t}}^{2}) \prod_{j=1}^{N_t}
f(B_{\tau_j}^{1}-B_{\tau_j}^{2})\right],
\end{equation}
(with the convention that on $\{N_t=0\}$, the product is defined to
be $1$), where $B^1=(B_t^1)_{t \geq 0}$ and $B^2=(B_t^2)_{t \geq 0}$
are independent $d$-dimensional Brownian motions starting from $x$,
respectively $y$, $N=(N_t)_{t \geq 0}$ is an independent Poisson
process with rate $1$ and points $\tau_1<\tau_2<\ldots$, and
$$w(t,x)=\int_{\bR^d} p_{t}(x-y)u_0(y)dy, \quad \mbox{where} \
p_t(x)=\frac{1}{(2\pi t)^{d/2}} e^{-|x|^2/(2t)}.$$ Note that the
representation (\ref{Dalang-representation}) does not rely on the
entire Brownian path, but only on its values at the (random) points
$\tau_1,\tau_2, \ldots, \tau_{N_t}$. This property has allowed the
authors of \cite{DMT08} to generalize the representation to a large
class of s.p.d.e.'s, including the wave equation.

To see where the idea for this representation comes from, we recall
briefly the salient points leading to (\ref{Dalang-representation}).
If the solution of (\ref{heat-Dalang}) exists, then it is unique and
admits the Wiener chaos expansion: (see e.g. Proposition 4.1 of
\cite{DMT08})
\begin{equation}
\label{Wiener-chaos-u-tx}
u_{t,x}=w(t,x)+\sum_{n=1}^{\infty}I_n(f_n(\cdot,t,x)),
\end{equation}
where $I_n(f_n(\cdot,t,x))=\int_{([0,1] \times \bR^d)^n}f_n(t_1,x_1,
\ldots, t_n,x_n) dF_{t_1,x_1} \ldots dF_{t_n,x_n}$, and $f_n \in
\cP^{\otimes n}$ is a symmetric function given by:
\begin{equation}
\label{defin-of-fn} f_{n}(t_{1},x_{1},\ldots,t_{n},x_{n},t,x)
=\frac{1}{n!} w(t_{\rho(1)},x_{\rho(1)}) \prod
_{j=1}^{n}p_{t_{\rho(j+1)}-t_{\rho(j)}}(x_{\rho(j+1)} -x_{\rho(j)}).
\end{equation}
Here, $\rho$ is a permutation of $\{1, \ldots, n\}$ such that
$t_{\rho(1)}<t_{\rho(2)}<\ldots <t_{\rho(n)}$, $t_{\rho(n+1)}=t$ and
$x_{\rho (n+1)}=x$. By the orthogonality of the terms in the series
(\ref{Wiener-chaos-u-tx}),
$$E[u_{t,x}u_{t,y}]=w(t,x)w(t,y)+\sum_{n=1}^{\infty}J_n(t,x,y),$$
where
\begin{eqnarray*}
J_n(t,x,y)&=&n! \langle f_n(\cdot,t,x),
f_n(\cdot,t,y)\rangle_{\cP^{\otimes n}}^2
=\int_{T_n(t)} F(t_1,\ldots, t_n) d{\bf t}, \quad \mbox{and} \\
F(t_1, \ldots,t_n) &=& \int_{\bR^{2nd}}
\prod_{j=1}^{n}p_{t_{j+1}-t_{j}}(x_{j+1}-x_{j})
\prod_{j=1}^{n}p_{t_{j+1}-t_{j}}(y_{j+1}-y_{j})\\
& & w(t_1,x_1)w(t_1,y_1) \prod_{j=1}^{n}f(x_j-y_j) d{\bf y} d{\bf x}
d{\bf t}.
\end{eqnarray*}
Here, we denote ${\bf t}=(t_1,\ldots,t_n), {\bf
x}=(x_1,\ldots,x_n)$, ${\bf y}=(y_1,\ldots,y_n)$ and $T_n(t)=\{
(t_1,\ldots,t_n);0<t_1<\ldots<t_n<t\}$.

A crucial observation of \cite{DMT08} is that, for any $F:[0,t]^{n}
\to \bR_{+}$ measurable,
\begin{equation}
\label{integral-over-symplex} \int_{T_n(t)} F(t_1,\ldots, t_n) d{\bf
t}=e^t E^N[F(t-\tau_n, \ldots, t-\tau_1)1_{\{N_t=n\}}].
\end{equation}
This is a key idea, which yields a probabilistic representation for
$J_n(t,x,y)$, based on the jump times of the Poisson process. This
idea is new in the literature, although it appeared implicitly in
the earlier works \cite{hersch74}, \cite{kac74}, \cite{pimsky91}.
Secondly, for any $0<t_1<\ldots<t_n<t$, $F(t-t_n, \ldots, t-t_1)$ is
represented as:
\begin{eqnarray}
\nonumber \lefteqn{F(t-t_n, \ldots, t-t_1) = } \\
\label{repr-of-F}
& & E^{B^1,B^2}\left[
w(t-t_n,B_{t_n}^1)w(t-t_n,B_{t_n}^2)\prod_{j=1}^{n}
f(B_{t_j}^{1}-B_{t_j}^{2})\right].
\end{eqnarray}
Relation (\ref{Dalang-representation}) follows from these
observations, using the independence between $N$ and $B^1,B^2$.

\vspace{3mm}

In this article, we generalize these ideas to the case of the
stochastic heat equation driven by a fractional-colored noise. More
precisely, we consider the following equation:
\begin{eqnarray}
\label{heat} \frac{\partial u}{\partial t} &=& \frac{1}{2} \Delta u
+ u \diamond \dot {W} , \quad t \in [0,1], x\in \bR^{d} \\
\nonumber u_{0,x} &=& u_0(x), \quad x \in \bR^d,
\end{eqnarray}
where $u_{0} \in C_b(\bR^d)$ is non-random, $\diamond$ denotes the
Wick product, and $\dot W$ is a Gaussian noise whose covariance is
formally given by:
$$E[\dot{W}_{t,x} \dot{W}_{s,y}]=\eta(t,s)f(x-y),$$
with $$\eta(t,s)=\alpha_H|t-s|^{2H-2}, \quad f(x)= \int_{\bR^{d}}
e^{-i \xi \cdot x} \mu (d\xi).$$ Here, $H \in (1/2,1)$,
$\alpha_H=H(2H-1)$ and $\mu$ is a tempered measure on $\bR^d$. (Note
that $\eta(t,s)$ is the covariance kernel of the fractional Brownian
motion.)

The noise $\dot{W}$ is defined rigourously as in
\cite{balan-tudor08a}, by considering a zero-mean Gaussian process
$W=\{W(h); h \in \cH \cP\}$ with covariance $$E(W(h)W(g))=\langle h,
g \rangle_{\cH \cP},$$ where $\cH \cP$ is the the completion of
$\{1_{[0,t] \times A}; t \in [0,1], A \in \cB_b(\bR^d) \}$ with
respect to the inner product  $\langle \cdot,\cdot \rangle_{\cH
\cP}$ given by:
$$\langle 1_{[0,t] \times A}, 1_{[0,s] \times B} \rangle_{\cH \cP}=
\int_0^t \int_0^s \int_{A} \int_{B} \eta(u,v)f(x-y)dy dxdv du.$$
(See also \cite{dalang99} for a martingale treatment of the case
$\eta(t,s)=\delta_0(t-s)$, which corresponds to $H=1/2$.)

As in \cite{balan-tudor08b}, the solution of equation (\ref{heat})
is interpreted in the mild sense, using the Skorohod integral with
respect to $W$. More precisely, an adapted square-integrable process
$u=\{u_{t,x}; (t,x) \in [0,1] \times \bR^{d}\}$ is a {\bf solution
to} (\ref{heat}) if for any $(t,x) \in [0,1] \times \bR^d$, the
process $\{p_{t-s}(x-y) u_{s,y} 1_{[0,t]}(s); (s,y) \in [0,1] \times
\bR^d\}$ is Skorohod integrable, and
$$u_{t,x}=p_t
u_0(x)+\int_{0}^{T}\int_{\bR^d}p_{t-s}(x-y)u_{s,y}1_{[0,t]}(s)
\delta W_{s,y}.$$

The existence of a solution $u$ (in the space of square-integrable
processes) depends on the roughness of the noise, introduced by $H$
and the kernel $f$. We mention briefly several cases which have been
studied in the literature. If $f=\delta_0$, the solution exists for
$d=1,2$ (see \cite{hu-nualart08}). If
$f(x)=\prod_{i=1}^{d}\alpha_{H_i}|x_i|^{2H_i-2}$, the solution
exists for $d< 2/(2H-1)+\sum_{i=1}^{d}H_i$ (see \cite{hu01}).
Recently, it was shown in \cite{balan-tudor08b} that if $f$ is the
Riesz kernel of order $\alpha$, or the Bessel kernel of order
$\alpha<d$, then the solution exists for $d \leq 2+\alpha$, whereas
if $f$ is the heat kernel or the Poisson kernel, then the solution
exists for any $d$.

\vspace{3mm}

Our main result establishes a probabilistic representation for the
second moment of the solution, similar to
(\ref{Dalang-representation}), and based on the planar Poisson
process.

\begin{theorem}
\label{main} Suppose that equation (\ref{heat}) has a solution
$u=\{u_{t,x}; (t,x) \in [0,1] \times \bR^d\}$ in $[0,1] \times
\bR^d$. Then, for any $(t,x) \in [0,1] \times \bR^d$ and $(s,y) \in
[0,1] \times \bR^d$,
$$E[u_{t,x}u_{s,y}]=w(t,x)w(s,y)+e^{ts} $$
$$\sum_{i_1, \ldots, i_n}E_{x,y}\left[w(t-\tau^*,B_{\tau^*}^{1})
w(s-\rho^*,B_{\rho^*}^2) \prod_{j=1}^{N_{t,s}}\eta (t-\tau_{i_j},
s-\rho_{i_j})f(B_{\tau_{i_j}}^{1}-B_{\rho_{i_j}}^{2})
I_{A_{i_1,\ldots,i_n}}\right],$$ where the sum is taken over all
distinct indices $i_1, \ldots,i_n$,
\begin{itemize}
\item $B^1=(B_t^1)_{t \in [0,1]}$ and $B^2=(B_t^2)_{t \in [0,1]}$ are
independent $d$-dimensional standard Brownian motions, starting from
$x$, respectively $y$,

\item $N=(N_{t,s})_{(t,s) \in [0,1]^2}$ is an independent planar
Poisson process with rate $1$ and points $\{P_i\}_{i \geq 1}$, with
$P_i=(\tau_i,\rho_i)$,

\item $A_{i_1, \ldots, i_n}(t,s)$ is the event that $N$ has points
$P_{i_1}, \ldots, P_{i_n}$ in $[0,t] \times [0,s]$,
$\tau^*=\max\{\tau_{i_1}, \ldots, \tau_{i_n}\}$ and
$\rho^*=\max\{\rho_{i_1}, \ldots,\rho_{i_n}\}$.
\end{itemize}
\end{theorem}

\section{Proof of Theorem \ref{main}}

We begin by recalling some basic facts about the planar Poisson
process. If $N=(N_{t,s})_{(t,s) \in [0,1]^2}$ is a $2$-parameter
process, and $R=(a,b] \times (c,d]$ is a rectangle in $[0,1]^2$, we
define $N_{R}=N_{a,c}+N_{b,d}-N_{a,d}-N_{b,c}$.

We say that $N=(N_{t,s})_{(t,s) \in [0,1]^2}$ is a {\bf planar
Poisson process} of rate $\lambda>0$, if it satisfies the following
conditions:

{\em (i)} $N$ vanishes on the axes, i.e. $N_{t,0}=N_{0,t}=0$ for all
$t \in [0,1]$.

{\em (ii)} $N_{R}$ has a Poisson distribution with mean $\lambda
|R|$, for any rectangle $R$;

{\em (iii)} $N_{R_1}, \ldots, N_{R_k}$ are independent, for any
disjoint rectangles $R_1, \ldots,R_k$.

\vspace{1mm}

The following construction of the planar Poisson process is
well-known (see e.g. \cite{AMSW83}). Let $X$ a Poisson random
variable with mean $\lambda$, and $\{P_i\}_{i \geq 1}$ be an
independent sequence of i.i.d. random vectors, uniformly distributed
on $[0,1]^2$. We denote $P_i=(\tau_i,\rho_i)$. For any $(t,s) \in
[0,1]^2$, define
$$N_{t,s}=\sum_{i=1}^{X}I\{\tau_i \leq t, \rho_{i} \leq s\}.$$
Then $N=(N_{t,s})_{(t,s) \in [0,1]^2}$ is a planar Poisson process
with rate $\lambda$.

For any $n \geq 1$ and for any distinct positive integers $i_1,
\ldots,i_n$, let $A_{i_1, \ldots,i_n}(t,s)$ be the event that $N$
has points $P_{i_1}, \ldots, P_{i_n}$ in $[0,t] \times [0,s]$. Then
\begin{equation}
\label{Nts-equal-n} \{N_{t,s}=n\}=\bigcup_{i_1, \ldots,i_n \geq 1 \
{\rm distinct}} A_{i_1, \ldots, i_n}(t,s).
\end{equation}

The following result is probably well-known. We include it for the
sake of completeness.

\begin{lemma}
The conditional distribution of  $$((t-\tau_{i_1},s-\rho_{i_1}),
\ldots, (t-\tau_{i_n},s-\rho_{i_n})) \quad \mbox{given} \quad
A_{i_1, \ldots,i_n}(t,s)$$ is uniform over $([0,t] \times [0,s])^n$.
\end{lemma}

\noindent {\bf Proof:} Let $I=\{i_1,\ldots,i_n\}$ and
$i^*=\max\{i_1, \ldots,i_n\}$. Then
$$A_{I}(t,s)=
\bigcup_{k \geq n, I \subset \{1, \ldots,k\}}  \{X=k\} \bigcap
\left(\bigcap_{i \in I}\{\tau_{i} \leq t, \rho_{i} \leq s\}\right)
\bigcap \left(\bigcap_{i \leq k,i \not \in I}\{\tau_i
>t \ {\rm or} \ \rho_i >s\} \right)$$
and $$P(A_I(t,s))=  \sum_{k \geq n \vee i^*}
e^{-\lambda}\frac{\lambda^k}{k!} (ts)^n(1-ts)^{k-n}.$$

For any Borel sets $\Gamma_1, \ldots, \Gamma_n$ in $[0,t] \times
[0,s]$, we have
$$P\left(\bigcap_{j=1}^{n}\{ (t-\tau_{i_j}, s-\rho_{i_j}) \in \Gamma_j\}
 \bigcap
A_{I}(t,s)\right)=\sum_{k \geq n \vee i^*} e^{-\lambda}
 \frac{\lambda^k}{k!}
\left(\prod_{j=1}^{n}|\Gamma_j|\right)(1-ts)^{k-n},$$ and hence
$$P\left(\bigcap_{j=1}^{n}\{(t-\tau_{i_j}, s-\rho_{i_j}) \in \Gamma_j\} \ | \
A_{I}(t,s) \right)=\frac{1}{(ts)^n} \prod_{j=1}^{n}|\Gamma_j|.$$
$\Box$

As a consequence, for any measurable function $F:([0,t] \times
[0,s])^n \to {\bf R}_{+}$,
\begin{equation}
\label{exp-F} E[F(t-\tau_{i_1},s-\rho_{i_1},
\ldots,t-\tau_{i_n},s-\rho_{i_n}) \ | \ A_{i_1, \ldots,i_n}(t,s)]=
\end{equation}
$$\frac{1}{(ts)^n} \int_{[0,t]^n} \int_{[0,s]^n} F(t_1,s_1, \ldots,
t_n,s_n) d{\bf s} d{\bf t},$$ where ${\bf t}=(t_1, \ldots,t_n)$ and
${\bf s}=(s_1, \ldots, s_n)$.

Suppose that $\lambda=1$. Then $(ts)^n=n! \ e^{ts}P(N_{t,s}=n)$ and
$$E[F(t-\tau_{i_1},s-\rho_{i_1}, \ldots,t-\tau_{i_n},s-\rho_{i_n}) \ | \
A_{i_1, \ldots,i_n}(t,s)]=$$ $$\frac{1}{n! e^{ts}P(N_{t,s}=n)}
\int_{[0,t]^n} \int_{[0,s]^n} F(t_1,s_1, \ldots, t_n,s_n) d{\bf s}
d{\bf t}.$$


Using (\ref{Nts-equal-n}) and (\ref{exp-F}), we obtain that for any
$F:([0,t] \times [0,s])^n \to {\bf R}_{+}$ measurable,
\begin{equation}
\label{integral-using-PPP} \int_{[0,t]^n} \int_{[0,s]^n} F(t_1,s_1,
\ldots, t_n,s_n) d{\bf s} d{\bf t} =
\end{equation}
$$ n! \ e^{ts} \sum_{i_1, \ldots, i_n \geq 1 \ {\rm distinct}} E^N [
F(t-\tau_{i_1},s-\rho_{i_1}, \ldots,t-\tau_{i_n},s-\rho_{i_n})
I_{A_{i_1, \ldots, i_n}(t,s)}].$$
Relation (\ref{integral-using-PPP}) is the analogue of
(\ref{integral-over-symplex}), needed in the fractional case.

\vspace{2mm}

Suppose now that equation (\ref{heat}) has a solution $u=\{u_{t,x};
(t,x) \in [0,1] \times \bR^d\}$. Then the solution is unique and
admits the Wiener chaos expansion (\ref{Wiener-chaos-u-tx}), where
$$I_n(f_n(\cdot,t,x))=\int_{([0,1] \times \bR^d)^n}f_n(t_1,x_1,
\ldots, t_n,x_n) dW_{t_1,x_1} \ldots dW_{t_n,x_n}$$ is the multiple
Wiener integral with respect to $W$, and $f_n \in \cH \cP^{\otimes
n}$ is the symmetric function given by (\ref{defin-of-fn}). (See
(7.4) of \cite{hu01}, or (4.4) of \cite{hu-nualart08}, or
Proposition 3.2 of \cite{balan-tudor08b}).

Using (\ref{Wiener-chaos-u-tx}), and the orthogonality of the Wiener
chaos spaces, we conclude that for any $(t,x) \in [0,1] \times
\bR^d$ and $(s,y) \in [0,1] \times \bR^d$,
\begin{eqnarray}
\nonumber E[u_{t,x}u_{s,y}]&=& E(u_{t,x})E(u_{s,y})+\sum_{n \geq
1}E[I_n(f_n(\cdot,t,x))
I_{n}(f_n(\cdot,s,y))] \\
\label{moment-utx-usy} &=& w(t,x)w(s,y)+ \sum_{n \geq 1}
\frac{1}{n!}\alpha_n(t,x,s,y),
\end{eqnarray}
 where $\alpha_n(t,x,s,y) := (n!)^2
\langle f_n(\cdot,t,x), f_n(\cdot,s,y) \rangle_{\cH \cP^{\otimes
n}}$ for any $n \geq 1$. Note that
\begin{equation}
\label{def-alpha-n} \alpha_n(t,x,s,y) = \int_{[0,t]^n}
\int_{[0,s]^n} \prod_{j=1}^{n}\eta(t_j,s_j) \langle G_{{\bf t};x},
G_{{\bf s};y} \rangle_{\cP(\bR^d)^{\otimes n}} d{\bf s} d{\bf t},
\end{equation}
where $\cP(\bR^d)$ is the the completion of $\{1_{A}; A \in
\cB_b(\bR^d) \}$ with respect to the inner product $\langle 1_{A},
1_{B} \rangle_{\cP(\bR^d)}= \int_{A} \int_{B} f(x-y)dy dx$,
\begin{eqnarray*}
G_{{\bf t};x}(x_1, \ldots x_n) & = & w(t_{\rho(1)}, x_{\rho(1)})
 \prod_{j=1}^{n}p_{t_{\rho(j+1)}-
t_{\rho}(j)}(x_{\rho(j+1)}-x_{\rho(j)}) \\
G_{{\bf s};y}(y_1, \ldots y_n) & = & w(s_{\sigma(1)}, y_{\sigma(1)})
\prod_{j=1}^{n} p_{s_{\sigma(j+1)}-
s_{\sigma}(j)}(y_{\sigma(j+1)}-y_{\sigma(j)}),
\end{eqnarray*}
the permutations $\rho$ and $\sigma$ are chosen such that:
$$0<t_{\rho(1)}<t_{\rho(2)}< \ldots < t_{\rho(n)} \quad \mbox{and}
\quad 0<s_{\sigma(1)}<s_{\sigma(2)}< \ldots < s_{\sigma(n)},$$
$t_{\rho(n+1)}=t$, $s_{\sigma(n+1)}=s$, $x_{\rho(n+1)}=x$ and
$y_{\sigma(n+1)}=y$.

\vspace{3mm}

Using (\ref{integral-using-PPP}) and (\ref{def-alpha-n}), we
conclude that $\alpha_n(t)$ admits the representation:
$$\alpha_{n}(t,x,s,y)=n! \ e^{ts}  \sum_{i_1, \ldots,i_n \geq 1 \
{\rm distinct}} E^{N}
\left[\prod_{j=1}^{n}\eta(t-\tau_{i_j},s-\rho_{i_j}) \right.$$
\begin{equation}
\label{alpha-n-tx-sy} \left. \langle G_{t-\tau_{i_1}, \ldots,
t-\tau_{i_n};x}, G_{s-\rho_{i_1}, \ldots, s-\rho_{i_n};y}
\rangle_{\cP(\bR^d)^{\otimes n}} I_{A_{i_1, \ldots,
i_n}(t,s)}\right].
\end{equation}


The next result gives a probabilistic representation of the inner
product
appearing in (\ref{alpha-n-tx-sy}). This result is the analogue of
(\ref{repr-of-F}), needed for the treatment of the fractional case.

\begin{lemma}
\label{representation-HP-product} For any $t_1, \ldots, t_n \in
[0,t]$, $s_1, \ldots, s_n \in [0,s]$, $x \in \bR^d$ and $y \in
\bR^d$,
\begin{eqnarray*}
\lefteqn{\langle G_{t-t_1, \ldots, t-t_n;x}, G_{s-s_1, \ldots,
s-s_n;y} \rangle_{\cP(\bR^d)^{\otimes n}}=} \\
& & E^{B^1,B^2} \left[w(t-t^*,B_{t^*}^{1}) w(s-s^*, B_{s^*}^{2})
\prod_{j=1}^{n} f(B_{t_j}^{1}-B_{s_j}^{2}) \right],
\end{eqnarray*} where
$t^*=\max\{t_1, \ldots,t_n \}$, $s^*=\max\{s_1, \ldots,s_n \}$, and
$B^1=(B_t^1)_{t \in [0,1]}$ and $B^{2}=(B_t^2)_{t \in [0,1]}$ are
independent $d$-dimensional Brownian motions, starting from $x$,
respectively $y$.
\end{lemma}

\noindent {\bf Proof:} Let $B^1$ and $B^2$ be independent
$d$-dimensional Brownian motions, starting from $0$. Note that
$$I:=\langle G_{t-t_1, \ldots, t-t_n;x}, G_{s-s_1, \ldots,
s-s_n;y} \rangle_{\cP(\bR^d)^{\otimes n}}=\int_{\bR^{2nd}}
w(t-t_{\rho(1)},x_{\rho(1)})w(s-s_{\sigma(1)},y_{\sigma(1)})$$
$$\prod_{j=1}^{n}f(x_j-y_j)p_{t_{\rho(j)}-t_{\rho(j+1)}}
(x_{\rho(j)}-x_{\rho(j+1)}) p_{s_{\sigma(j)}- s_{\sigma(j+1)}}
(y_{\sigma(j)}-y_{\sigma(j+1)})d{\bf y}d{\bf x},$$ where the
permutations $\rho$ and $\sigma$ of $\{1,\ldots,n\}$ are chosen such
that
$$0<t_{\rho(n)}< \ldots < t_{\rho(1)}<t \quad \mbox{and}
\quad 0<s_{\sigma(n)}< \ldots < s_{\sigma(1)}<s,$$
$t_{\rho(n+1)}=s_{\sigma(n+1)}=0$, $x_{\rho(n+1)}=x$,
$y_{\sigma(n+1)}=y$. We use the change of variables
$$x_{\rho(j)}-x_{\rho(j+1)}=z_{n+1-j}, \quad
y_{\sigma(j)}-y_{\sigma(j+1)}=w_{n+1-j}, \quad  j=1, \ldots,n.$$
Note that $x_{\rho(j)}=x+\sum_{k=1}^{n+1-j}z_k$, i.e.
$x_{j}=x+\sum_{k=1}^{n+1-\rho^{-1}(j)}z_k$. We get:
$$I=\int_{\bR^{2nd}}
w(t-t_{\rho(1)},x+\sum_{k=1}^{n}z_k)w(s-s_{\sigma(1)},y+\sum_{k=1}^{n}w_k)$$
$$\prod_{j=1}^{n}f\left((x+
\sum_{k=1}^{n+1-\rho^{-1}(j)}z_k)-(y+\sum_{k=1}^{n+1-\sigma^{-1}(j)}w_k)\right)$$
$$\prod_{j=1}^{n}p_{t_{\rho(j)}-t_{\rho(j+1)}}(z_{n+1-j}) \prod_{j=1}^{n}
p_{s_{\sigma(j)}- s_{\sigma(j+1)}} (w_{n+1-j}) d{\bf w}d{\bf z}.$$

We now use the fact that $p_{t-s}(x)dx$ is the density of the
increment of a $d$-dimensional Brownian motion over the interval
$(s,t]$, these increments over disjoint intervals are independent,
and the Brownian motions $B^1,B^2$ are independent. Therefore, we
replace the integral over $\bR^{2nd}$ by the expectation
$E^{B^1,B^2}$, and the variables $z_k,w_k$ by
$B_{t_{\rho(n+1-k)}}^{1}- B_{t_{\rho(n+1-k+1)}}^{1}$, respectively
$B_{s_{\sigma(n+1-k)}}^{2}- B_{s_{\sigma(n+1-k+1)}}^{2}$, for all
$k=1,\ldots,n$. Note that, for any $m=1,\ldots,n$
$$\sum_{k=1}^{m}(B_{t_{\rho(n+1-k)}}^{1}-
B_{t_{\rho(n+1-k+1)}}^{1})=B_{t_{\rho(n+1-m)}}^{1}$$
$$\sum_{k=1}^{m}(B_{s_{\sigma(n+1-k)}}^{2}-
B_{s_{\sigma(n+1-k+1)}}^{2})=B_{s_{\sigma(n+1-m)}}^{2}.$$ Hence,
$$I=E^{B^1,B^2}
[w(t-t_{\rho(1)},x+B_{t_{\rho(1)}}^{1})w(s-s_{\sigma(1)},y+
B_{s_{\sigma(1)}}^{2})\prod_{j=1}^{n}f((x+B_{t_j}^{1})-(y+B_{s_j}^{2}))].$$
$\Box$

\noindent {\bf Conclusion of the Proof of Theorem \ref{main}:} Since
equation (\ref{heat}) has a solution, this solution is unique and
admits the Wiener chaos expansion (\ref{Wiener-chaos-u-tx}). The
result follows from (\ref{moment-utx-usy}), (\ref{alpha-n-tx-sy})
and Lemma \ref{representation-HP-product}. $\Box$

\end{document}